\documentclass[11pt,leqno]{article} 
\usepackage{graphics}
\newtheorem{thm}{Theorem}[section]
\newtheorem{lma}{Lemma}[section]

\newcommand{\beqa}{\begin{eqnarray}}
\newcommand{\eeqa}{\end{eqnarray}}

\newcommand{\pf}{\noindent {\bf Proof:} $\s$ }
\newcommand{\epf}{ \hfill$\diamondsuit$ \medskip}

\newcommand{\beq}{\begin{equation}}
\newcommand{\eeq}{\end{equation}}
\newcommand{\lbl}{\label}
\newcommand{\s}{\; \;}

\newcommand{\ep}{\epsilon}

\newcommand{\la}{\lambda}

\newcommand{\ra}{\rightarrow}
\newcommand{\al}{\alpha}

\title{On the perturbed Gelfand equation from combustion theory}

\author{
Philip Korman \thanks{Supported in part by the Taft faculty grant at the University of Cincinnati}  \\ 
Department of Mathematical Sciences \\ 
University of Cincinnati \\ 
Cincinnati Ohio 45221-0025 \\
\\
Yi Li\\
Department of Mathematics \\
California State University, Northridge \\
Northridge, CA 91330-8313\\
\\
Tiancheng Ouyang \\
Department of Mathematics \\
Brigham Young University \\
Provo, Utah  84602
}
\begin{document}
\maketitle
\begin{abstract} 
For the  perturbed Gelfand's equation on the unit ball in two dimensions, 
Y. Du and Y. Lou \cite{DL} proved that the curve of positive solutions is exactly $S$-shaped, for sufficiently small values of the secondary parameter. We present a simplified proof and some extensions. This problem is prominent in combustion theory, see e.g., the book of J. Bebernes and D. Eberly \cite{BE}.
 \end{abstract}

\begin{flushleft}
\end{flushleft}

\begin{flushleft}
\end{flushleft}

\section{Introduction}
\setcounter{equation}{0}
The following Dirichlet problem for the  perturbed Gelfand's equation is prominent in combustion theory
\beq
\lbl{100}
\Delta u+\la e^{\frac{u}{1+\ep u}} =0 \,, \s \mbox{for $|x|<1$}\,, \s  u=0 \s\s  \mbox{when $|x|=1$} \,,
\end{equation}
see e.g., J. Bebernes and D. Eberly \cite{BE}. Here $\lambda$ and $\epsilon$ are positive parameters, and we think of $\lambda$ as the primary parameter, while $\epsilon$ is the secondary, or ``evolution parameter". By the maximum principle, the solution of (\ref{100}) is positive, and then by the classical theorem of B. Gidas, W.-M. Ni and L. Nirenberg \cite{GNN} it is radially symmetric, i.e., $u=u(r)$, with $r=|x|$, and it satisfies
\[
u'' +\frac{n-1}{r}u'+\la e^{\frac{u}{1+\ep u}} =0 \,, \s \mbox{$0<r<1$}\,, \s  \s u'(0)=u(1)=0 \,.
\]
For the perturbed Gelfand's equation on a unit ball in two dimensions
\begin{equation}
\label{testi10}
u'' +\frac{1}{r}u'+\la e^{\frac{u}{1+\ep u}} =0 \,, \s \mbox{$0<r<1$}\,, \s  \s u'(0)=u(1)=0 \,,
\end{equation}
Y. Du and Y. Lou \cite{DL}, building on the earlier results of P. Korman, Y. Li and T. Ouyang \cite{KLO}, \cite{KLO2}, proved the following theorem, thus settling a long-standing conjecture of S.V. Parter  \cite{P}.
\begin{thm}\label{thm:*}
For $\ep$ sufficiently small the solution set of (\ref{testi10}) is exactly $S$-shaped.  Moreover, at any $\lambda$ where either two or three solutions occur, these solutions are strictly ordered. (See Figure $1$.)
\end{thm}

Their proof was rather involved, and it was relying on some previous results of E.N. Dancer \cite{D}. The purpose of this note is to give a simpler and self-contained proof. We also observe that some more general results can be obtained without too much extra effort. While simplifying the proof in \cite{DL}, we retain several of the crucial steps from that paper: the change of variables, Lemma \ref{lma:test10} and Theorem \ref{thm:3.1}, although we generalize or simplify these results. Our new tool involves showing that the turning points are {\em non-degenerate}, so that they persist when the  secondary parameter $\epsilon$ is varied. We also observe that computer assisted validation of bifurcation diagrams is possible for $\epsilon$ not being small.

\begin{figure}
\begin{center}
\scalebox{0.8}{\includegraphics{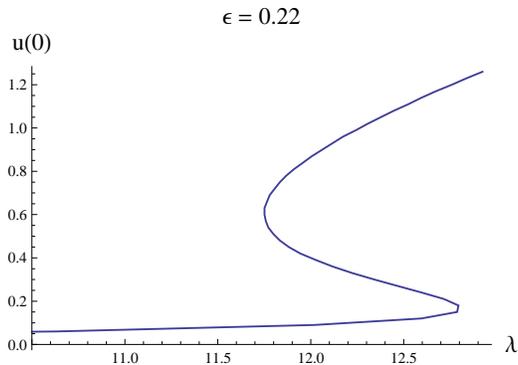}}
\end{center}
\caption{ An $S$-shaped solution curve when $\epsilon=0.22$}
\end{figure}
\medskip

\medskip

When $\epsilon \geq 0.25$, the solution curve of (\ref{testi10}) is monotone. Indeed, in that case the function $f(u)=e^{\frac{u}{1+\ep u}}$ satisfies $uf'(u)<f(u)$ for all $u>0$, except for $u=4$ when $\epsilon = 0.25$, so that the corresponding linearized problem has only the trivial solution by the Sturm's comparison theorem, and so the implicit function theorem applies. It is natural to conjecture that there is a critical $\epsilon _0>0$, so that for $\epsilon \geq \epsilon _0>0$ the solution curve  is monotone, while for $\epsilon < \epsilon _0>0$ it is exactly  $S$-shaped. For the one-dimensional case the same statement is  known as S.-H. Wang's conjecture \cite{W1}, for which we gave a computer assisted proof in \cite{KLO3}. Our numerical computations, for the two-dimensional case  (\ref{testi10}), show that remarkably $\epsilon _0$ is  close to $0.25$, see Figures $1$ and $2$. (We used the shoot-and-scale method, described in detail in \cite{K2}, which we implemented using {\em Mathematica}.) Previous contributions to the $n=1$ case included K.J. Brown et al \cite{Bis}, S.P. Hastings and J.B. McLeod  \cite{H} (who proved the above theorem in one dimension), S.-H. Wang \cite{W1},  S.-H. Wang  and F.P. Lee  \cite{W2}, P. Korman and Y. Li \cite{KL1}.
\medskip

When $n>2$, this theorem does  not hold. Indeed, if $\ep=0$ and $3 \leq n \leq 9$, by the classical result of  D.D. Joseph and  T.S. Lundgren  
\cite{JL} the solution curve of the problem (\ref{testi10}) makes infinitely many turns. It is natural to expect that these turns will persist for small $\ep >0$, which is indeed the case, as shown by E.N. Dancer \cite{D}.

\begin{figure}
\begin{center}
\scalebox{0.8}{\includegraphics{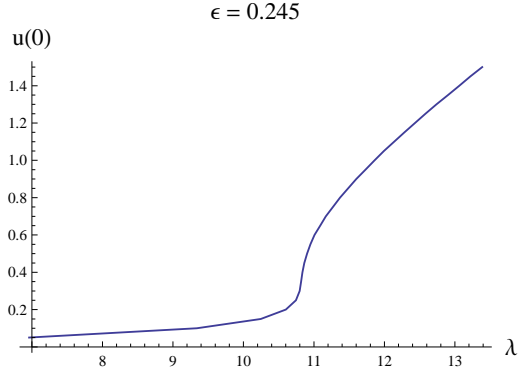}}
\end{center}
\caption{ A monotone  solution curve when $\epsilon=0.245$}
\end{figure}
\medskip

\section{Non-degenerate critical points}
\setcounter{equation}{0}

We consider positive solutions of the Dirichlet problem
\begin{equation}
\label{dul1}
\Delta u+\lambda f(u, \epsilon)=0 \s\s \mbox{for $x \in D$}, \s u=0 \s \mbox{on $\partial D$} \,,
\end{equation}
depending on positive parameters $\la$ and $\ep$, on a general domain $D \subset R^n$.
We are
interested in the {\em critical points} of (\ref{dul1}), i.e., the solution triples $(\la, u, \ep)$, for which
the corresponding linearized problem
\begin{equation}
\label{dul2}
\Delta w+\lambda f_u(u, \epsilon)w=0 \s\s \mbox{for $x \in D$}, \s w=0 \s \mbox{on $\partial D$} 
\end{equation}
has a nontrivial solution $w(x)$. We  wish to continue the critical points, when the secondary parameter $\ep$ is
varied. The following lemma was first proved by E.N. Dancer \cite{D} (see also J. Shi
\cite{S2}). We present a simpler proof for completeness.
\begin{lma}\label{lma:test0}
 Let $(\la _0,u_0, \ep _0)$ be a critical point of (\ref{dul1}). Assume that the null-space of (\ref{dul2}) is one-dimensional, spanned by some  $w_0(x) \in W^{2,2}(D) \cap W^{1,2}_0(D)$.
Assume also that
\begin{equation}
\label{dul3a}
\int_D f(u_0(x))w_0(x) \, dx \ne 0 \,,
\end{equation}
\begin{equation}
\label{dul3}
\int_D f''(u_0(x))w_0^3(x) \, dx \ne 0 \,.
\end{equation}
Then there is a unique critical point $(\la (\ep),u(\ep),\ep)$ near $(\la _0,u_0, \ep _0)$. Moreover, these are the only  critical points in some neighborhood of  $(\la _0,u_0, \ep _0)$.
\end{lma}

\pf
We can normalize the solution of (\ref{dul2}), so that
\begin{equation}
\label{dul4}
\frac12 \int_D w^2(x) \, dx =1 \,.
\end{equation}
The equations (\ref{dul1}), (\ref{dul2}) and (\ref{dul4}) give us three equations to find $u$, $w$ and $\la$ as a function of $\ep$. We show that the implicit function theorem applies. Indeed, we define
a map 
$
H(u,w,\la,\ep): \left(W^{2,2}(D) \cap W^{1,2}_0(D) \right) \times  \left(W^{2,2}(D) \cap W^{1,2}_0(D) \right) \times R \times R \rightarrow L^2(D) \times L^2(D) \times R
$
as a vector whose entries are the left hand sides of the above equations:
\[
H(u,w,\la,\ep)=\left[
\begin{array}{r}
\Delta u+\lambda f(u, \epsilon) \\
\Delta w+\lambda f_u(u, \epsilon)w \\
\frac12 \int_D w^2(x) \, dx
\end{array}
\right] \,.
\]
The linearized operator with respect to the first three variables at a point $(u_0,w_0,\la_0,\ep _0)$ is
\[
H'_{(u,w,\la)}(u_0,w_0,\la_0,\ep _0)\left[
\begin{array}{r}
 v \\
\theta \\
\tau
\end{array}
\right]
\]
\[
=\left[
\begin{array}{c}
\Delta v+\lambda _0 f_u(u_0, \epsilon_0)v +\tau f(u_0, \epsilon_0)\\
\Delta \theta+\lambda _0 f_u(u_0, \epsilon_0) \theta+\lambda _0 f_{uu}(u_0, \epsilon_0)w_0v +\tau f_u(u_0, \epsilon_0) w_0\\
 \int_D w_0 \theta \, dx
\end{array}
\right] \,.
\]

We need to show that this operator is both injective and surjective. To see that it is
injective, we need to show that the system
\begin{eqnarray}
\label{dul5}
& \Delta v+\lambda _0 f_u(u_0, \epsilon_0)v +\tau f(u_0, \epsilon_0)=0 \\ \nonumber
& \Delta \theta+\lambda _0 f_u(u_0, \epsilon_0) \theta+\lambda _0 f_{uu}(u_0, \epsilon_0)w_0v +\tau f_u(u_0, \epsilon_0) w_0=0 \\ \nonumber
&  \int_D w_0 \theta \, dx=0
\end{eqnarray}
has only the trivial solution $(v, \theta, \tau) = (0,0,0)$. The first equation in (\ref{dul5}) can be
regarded as a linear equation for $v$, with its kernel spanned by $w_0$, and the right
hand side equal to $-\tau f(u_0, \epsilon_0)$. Since by (\ref{dul3a}), $f(u_0, \epsilon_0)$ is not orthogonal
to the kernel, it follows that the first equation is solvable only if $\tau = 0$. We then
have $v = kw_0$, with a constant $k$. We now regard the second equation in (\ref{dul5})
as a linear equation for $\theta$ with the same kernel spanned by $w_0$, and the right hand side equal to
$-k \la _0 f_u(u_0, \epsilon_0) w^2_0$. By our condition (\ref{dul3}), the second equation is solvable only if
$k = 0$. We then have $\theta=lw_0$, with a constant $l$. From the third equation in (\ref{dul5})
we conclude that $l=0$, completing the proof of injectivity.
\medskip

Turning to the surjectivity, we need to show that for any $L^2(D)$ functions $a(x)$ and
$b(x)$, and for any constant $c$ the problem
\begin{eqnarray}
\label{dul6}
& \Delta v+\lambda _0 f_u(u_0, \epsilon_0)v +\tau f(u_0, \epsilon_0)=a(x) \\ \nonumber
& \Delta \theta+\lambda _0 f_u(u_0, \epsilon_0) \theta+\lambda _0 f_{uu}(u_0, \epsilon_0)w_0v +\tau f_u(u_0, \epsilon_0) w_0=b(x) \\ \nonumber
&  \int_D w_0 \theta \, dx=c
\end{eqnarray}
is solvable. Proceeding similarly to the above, we regard the first equation in (\ref{dul6})
as a linear equation for $v$ with the right hand side equal to $a(x) - \tau f(u_0, \epsilon_0)$. By
 by (\ref{dul3a}),  we can choose $ \tau$, so that this function is orthogonal to $w_0$.  Then the
first equation in (\ref{dul6}) has infinitely many solutions of the form $v= \bar v  + kw_0$, where $\bar v$
is any fixed solution, and $k$ is any constant. We now turn to the second equation in
(\ref{dul6}), where $\tau$ has been just fixed above. In view of our condition (\ref{dul3}), we can fix
$k$ so that $b(x)-\lambda _0 f_{uu}(u_0, \epsilon_0)w_0v -\tau f_u(u_0, \epsilon_0) w_0$ is orthogonal to $w_0(x)$, and hence the second equation is solvable. We then have $\theta = \bar \theta + l w_0$, where $\bar \theta$ is a fixed
solution and $l$ is any number. Finally, from the third equation in (\ref{dul6}) we uniquely
determine $l$.
\epf

This lemma shows that the critical points continue on a smooth curve, when the secondary parameter $\epsilon$ varies. We call such critical points {\em non-degenerate}.

\section{Positivity for the linearized problem}
\setcounter{equation}{0}

Let $u(r)$ be a positive solution of  the Dirichlet problem
\begin{equation}
\label{test0}
u'' +\frac{n-1}{r}u'+f(u) =0 \,, \s \mbox{$r>0$}\,, \s  \s u'(0)=u(1)=0 \,.
\end{equation}
We wish to show that any non-trivial solution of the linearized problem
\begin{equation}
\label{test1}
L[w] \equiv w'' +\frac{n-1}{r}w'+f'(u)w =0 \,, \s \mbox{$r>0$}\,, \s  \s w'(0)=w(1)=0
\end{equation}
does not vanish on $(0,1)$.

\begin{lma}\label{lma:test1}
Assume that there is a function $z(r) \in C^2(0,1)$ (a ``test function"), such that for some $\xi \in [0,1]$
\begin{equation}
\label{test2}
z>0 \,,\s \mbox{and} \s L[z]=z'' +\frac{n-1}{r}z'+f'(u)z < 0 \s \mbox{on $(0,\xi)$} \,,
\end{equation}
\begin{equation}
\label{test3}
z<0 \,,\s \mbox{and} \s L[z]=z'' +\frac{n-1}{r}z'+f'(u)z > 0 \s \mbox{on $(\xi,1)$} \,.
\end{equation}
Then $w(r)$  does not vanish on $[0,1)$, i.e., we may assume that $w(r)>0$.
\end{lma}

\pf
Without loss of generality we may assume that $w(0)>0$. We claim that $w(r)$ cannot vanish on $(0,\xi]$. Assuming the contrary, we can find $\xi _0 \in (0,\xi]$, so that $w(r)>0$ on $(0,\xi_0)$ and $w(\xi_0)=0$, 
$w'(\xi_0)<0$. Combining the equations in (\ref{test1}) and (\ref{test2}), we get
\[
\left[r^{n-1} \left( w'z-wz' \right) \right]'>0 \,.
\]
Integration over $(0,\xi_0)$ gives
\[
\xi_0^{n-1} w'(\xi_0) z(\xi_0)>0 \,,
\]
but $w'(\xi _0)<0$ and $z(\xi _0) \geq 0$, which is a contradiction.
\medskip

We show similarly that $w(r)$ cannot vanish on $[\xi,1)$. 
\epf

In a nutshell, we showed that $z(r)$ oscillates faster than $w(r)$ on both of the intervals $(0,\xi]$ and $[\xi,1)$, and hence $w(r)$ cannot vanish on either of the intervals.  Observe that both the cases $\xi=0$ and $\xi=1$ are allowed. Only rarely can one use this lemma directly, but rather the idea of its proof is used.
\medskip

We now present a generalization of a result of Y. Du and Y. Lou \cite{DL}.

\begin{thm}\label{thm:3.1}
For the problem (\ref{test0}) assume that $n=2$, the function $f(u) \in C^2(\bar R_+)$ satisfies $f(u)>0$, $f'(u)>0$ for all $u>0$,  and it is log-concave, i.e., 
\beq
\lbl{test2a}
f''(u)f(u)-{f'}^2(u)<0 \,, \;\; \mbox{for all $u>0$} \,.
\eeq
Let $u(r)$ be a positive solution of  (\ref{test0}).
Then  any non-trivial solution of the corresponding linearized problem (\ref{test1}) may be assumed to satisfy $w(r)>0$ on $[0,1)$.
\end{thm}

\noindent {\bf Proof:} $\s$ 
Write $f(u)=e^{h(u)}$, with $h'(u)>0$, and $h''(u)<0$ for all $u>0$, by (\ref{test2a}). Consider $z(r)=ru'(r) +\alpha$, with the constant $\alpha>0$ to be specified.  We have (using (\ref{test0}) with $n=2$)
\[
z'(r)=ru''(r)+u'(r)=-rf(u(r))<0 \,,
\]
so that $z(r)$ can vanish at most once on $(0,1)$. Compute
\[
L[z]=-2f(u)+\alpha f'(u)=e^{h(u)} \left[-2+\alpha h'(u) \right] \,.
\]
The function $g(u) \equiv -2+\alpha h'(u)$ satisfies
\[
\frac{d}{dr} g(u)=\alpha h''(u)u'(r)>0 \,,
\]
and hence $L[z]$ can change sign at most once on $(0,1)$ (from negative to positive).
\medskip

Without loss of generality we may assume that $w(0)>0$, and let us suppose that $w(\xi_0)=0$ at some $\xi_0 \in (0,1)$. Define $\alpha _1>0$ so that $z(r)=ru'(r) +\alpha_1$ also vanishes at $\xi_0$, and $z(r)>0$ on $(0,\xi_0)$. For $\alpha>0$ small, the function $g(u(r))=-2+\alpha  h'(u(r))$ is negative, while for larger $\alpha$'s it changes sign. Let $\alpha _2$ be the supremum of $\alpha$'s for which $g(u(r)) < 0 $ for all $r \in (0,1)$, so that for $\alpha> \alpha _2$, $g(u(r))$ changes sign (exactly once) on $(0,1)$.
\medskip

\noindent
Case 1. $\alpha _2 \geq \alpha _1$. Fix $z(r)=ru'(r) +\alpha_1$. Then on $(0,\xi_0)$, $z>0$ and $L[z]<0$, and hence $w(r)$ cannot vanish on $(0,\xi_0]$, a contradiction.
\medskip

\noindent
Case 2. $\alpha _2 < \alpha _1$. Then there is a point $\xi _1 \in (0,1)$, such that $g(u(r))=-2+\alpha _1  h'(u(r))<0$  on $(0,\xi_1)$ and $g(u(r))=-2+\alpha _1  h'(u(r))>0$  on $(\xi_1,1)$.  
\medskip

\noindent
Sub-case (a). $\xi _0 \leq \xi _1$. Again, fix $z(r)=ru'(r) +\alpha_1$. Then on $(0,\xi_0)$, $z>0$ and $L[z]<0$, and hence $w(r)$ cannot vanish on $(0,\xi_0]$, a contradiction.
\medskip

\noindent
Sub-case (b). $\xi _1 < \xi _0$. Consider $z(r)=ru'(r) +\alpha$, with $\alpha<\al_1$. As we decrease $\alpha$ from $\alpha _1$, the root of the decreasing function $z(r)$ moves to the left (and the root is at $r=0$, when $\alpha=0$), while the root of the increasing function $g(u(r))=-2+\alpha   h'(u(r))$ moves to the right. At some $
\bar \alpha>0$, these roots intersect at some $\bar \xi \in (\xi _1, \xi _0)$. Fix $z(r)=ru'(r) +\bar \alpha$. Then on $(0,\bar \xi)$, we have $z>0$ and $L[z]<0$,
and we have $z<0$ and $L[z]>0$ on $(\bar \xi,1)$. By Lemma \ref{lma:test1}, $w(r)$ cannot vanish (or a contradiction is achieved on $(\bar \xi,1)$).
\hfill$\diamondsuit$ \medskip

\noindent
{\bf Example}. $\displaystyle f(u)=u^p+u^q$, with positive constants $p \ne q$. A direct computation shows that $ f(u)$ is log-concave for all $u>0$, and the theorem applies, if and only if
\[
(p-q)^2-2(p+q)+1<0 \,.
\]

\noindent
{\bf Example}. $\displaystyle f(u)=e^{-\frac{1}{u+a}}$, with a constant $a \geq 0$. This function is log-concave for all $u>0$, and the theorem applies.
\medskip

This result does not hold   in case $n>2$. 

\begin{thm}\lbl{thm:du20}
In addition to the conditions of Theorem \ref{thm:3.1}, assume that
\[
uf'(u)> f(u)\,, \;\; \mbox{for all $u>0$} \,.
\]
Then any positive solution of (\ref{test0}) is non-singular, i.e., (\ref{test1}) has only the trivial solution.
\end{thm}

\pf
By the Sturm comparison theorem, $w(r)$ must vanish on $(0,1)$, in contradiction with the  Theorem \ref{thm:3.1}.
\epf

\section{The limiting problem}
\setcounter{equation}{0}

We consider now the solutions of the perturbed Gelfand's problem
\begin{equation}
\label{test10}
u'' +\frac{1}{r}u'+\la e^{\frac{u}{1+\ep u}} =0 \,, \s \mbox{$0<r<1$}\,, \s  \s u'(0)=u(1)=0 \,,
\end{equation}
on two-dimensional unit ball, with  positive parameters $\lambda$ and $\epsilon$. 

Following Y. Du  and Y. Lou \cite{DL}, we set $w(r) =\ep ^2 u(r)$, and $\mu=\la \ep ^2 e^{1/\epsilon}$, converting this problem into
\begin{equation}
\label{test11}
w'' +\frac{1}{r}w'+ \mu e^{-\frac{1}{\ep + w}} =0 \,, \s \mbox{$0<r<1$}\,, \s  \s w'(0)=w(1)=0 \,.
\end{equation}
The limiting problem at $\ep=0$ is 
\begin{equation}
\label{test12}
v'' +\frac{1}{r}v'+ \eta e^{-\frac{1}{ v}} =0 \,, \s \mbox{$0<r<1$}\,, \s  \s v'(0)=v(1)=0 \,,
\end{equation}
where we changed the names of the variables for future reference. The exact multiplicity result for (\ref{test12}) will follow from the following theorem, which is not much harder to prove than the special case of the problem (\ref{test12}). Except for the last statement, it corresponds to Theorem 2  in Y. Du  and Y. Lou \cite{DL}.
\begin{thm}\label{thm:test10}
Consider the problem
\begin{equation}
\label{test14}
v'' +\frac{1}{r}v'+ \eta f(v) =0 \,, \s \mbox{$r>0$}\,, \s  \s v'(0)=v(1)=0 \,.
\end{equation}
Assume that the function $f(v) \in C^2(\bar R_+)$ satisfies $f(0)=f'(0)=0$, $f(v)>0$ for $v>0$. Assume also that $f(v)$ is log-concave and convex-concave for all $v>0$, so that (\ref{test2a}) holds, and $f''(v)>0$ on $(0,\beta)$ and $f''(v)<0$ on $(\beta ,\infty)$, for some $\beta >0$.  Assume that $\lim _{v \rightarrow \infty} f(v)=f_0>0$. Assume finally that 
\begin{equation}
\label{test15}
vf'(v)>f(v) \;\; \mbox{on $(0,\beta)$} \,.
\end{equation}
Then there is a critical $\eta _0$, such that for $\eta <\eta _0$ the problem (\ref{test14}) has no positive solutions, it has exactly one positive solution at $\eta =\eta _0$, and there are exactly two positive solutions for $\eta >\eta _0$. Moreover, all solutions lie on a single smooth solution curve, which for $\eta >\eta _0$ has two branches, denoted by $v^{-}(r, \eta) < v^{+}(r, \eta)$, with $v^{+}(r, \eta)$ strictly monotone increasing in $\eta$, and $\lim _{\eta \rightarrow \infty} v^{+}(r, \eta)=\infty$ for all $r \in [0,1)$. For the lower branch, $\lim _{\eta \rightarrow \infty} v^{-}(0, \eta)=0$. Denote $v_0=v^{-}(r, \eta _0) =v^{+}(r, \eta _0)$. The turning point $(\eta _0, v_0)$ is the only critical point on the solution curve, and moreover $(\eta _0, v_0)$ is a non-degenerate  critical point.
\end{thm}

\pf
We follow  \cite{DL} to prove that   the problem (\ref{test14}) has a positive solution for $\eta$ large. Let $\phi \in C_0^{\infty}(B)$, $\phi \geq 0$, $\max _B \phi >0$ ($B$ is the unit ball in $R^n$). Let $\underline{v}>0$ be the unique solution of $\Delta v+\phi=0$ in $B$, $v=0$ on $\partial B$, and denote by $\overline{v}$ be the unique solution of $\Delta v+\eta f_0=0$ in $B$, $v=0$ on $\partial B$. Then $\underline{v}< \overline{v}$, for $\eta$ large, and they form a lower-upper solution pair.

\medskip
We now continue this solution for decreasing $\eta$. Since $f'(0)=0$, the solution curve cannot enter the point $(0,0)$ in the $(\la,u(0))$ plane, nor can a bifurcation from zero occur at some $\lambda >0$ (just multiply the PDE version of  (\ref{test14}) by $v$ and integrate). Hence, the solution curve must turn to the right. By the condition (\ref{test15}), and Theorem \ref{thm:du20}, any turns must occur in the region where $v(0)>\beta$. By a result of P. Korman, Y.   Li  and  T. Ouyang \cite{KLO},  only turns to the right are possible in that region. Hence, there is only one turn to the right. In particular, the proof in \cite{KLO} showed that the condition (\ref{dul3}) holds, and so the turning point is non-degenerate. Monotonicity of the upper branch, $v^{+}(r, \eta)$, is proved as in \cite{KLO}. (Near the turning point the monotonicity follows by the Crandall-Rabinowitz theorem. Then one uses the maximum principle to show that $v^{+}_{\eta}(r, \eta)>0$, for all $\eta >\eta _0$.)
\epf

In particular, this theorem applies to the problem (\ref{test12}) (we define $e^{-\frac{1}{ v}}$ to be zero at $v=0$). 
\medskip

The following important lemma was proved first in Y. Du  and Y. Lou \cite{DL}. We present its proof for completeness.

\begin{lma}\label{lma:test10}
Fix $\ep <v_0(0)$. Consider the solutions of (\ref{test11}), with $w(0)>v_0(0)-\ep$. Then $\mu$ is an increasing function of $w(0)$, i.e., the solution curve of  (\ref{test11}) travels to the right (northeast) in the $(\mu ,w(0))$ plane.
\end{lma}

\pf
Denote by $v(r,\alpha)$  the solution of (\ref{test12}) satisfying $ v(0,\alpha)=\alpha$. Assume that $\beta>\alpha>v_0(0)$, and $v(r,\beta)$ is the solution  of (\ref{test12}) with $ v(0,\beta)=\beta$. Recall that the solutions of (\ref{test12}) are uniquely identified by a global parameter $v(0)$, so that we have solution pairs $(\eta(\alpha),v(r,\alpha))$ and $(\eta(\beta),v(r,\beta))$, with $\eta(\beta)>\eta(\alpha)$ and $v(r,\beta)>v(r,\alpha)$ for all $r \in [0,1)$. Denote by $a=a(\alpha)$ the point where $v(a,\alpha)=\ep$, and set $w(r)=v(r,\alpha)-\epsilon$. Then $w(r)>0$ on $[0,a)$, and it satisfies
\[
w'' +\frac{1}{r}w'+\eta(\alpha)e^{-\frac{1}{\ep + w}} =0 \,, \s \mbox{$0<r<a$}\,, \s  \s w'(0)=w(a)=0 \,.
\]
Scaling $r=at$, we see that $w(t)$ satisfies
\[
w'' +\frac{1}{t}w'+\eta(\alpha)a^2(\alpha)e^{-\frac{1}{\ep + w}} =0 \,, \s \mbox{$0<r<1$}\,, \s  \s w'(0)=w(1)=0 \,.
\]
We see that $w(t)$ is the solution of (\ref{test11}) with the maximum value of $\al-\ep$ and the corresponding parameter value $\mu(\al-\ep)=\eta(\alpha)a^2(\alpha)$. Similarly, we identify $w(r)=v(r,\beta)-\epsilon$ as the solution of (\ref{test11}) with the maximum value of $\beta-\ep$ and the corresponding parameter value $\mu(\beta-\ep)=\eta(\beta)a^2(\beta)$. Observe that $a(\beta)>a(\alpha)$, because $v(r,\beta)>v(r,\alpha)$, and then 
\[
\mu(\beta-\ep)=\eta(\beta)a^2(\beta)>\eta(\alpha)a^2(\alpha)=\mu(\alpha-\ep) \,,
\]
and the proof follows.
\epf

\section{Proof of the Theorem \ref{thm:*}}
\setcounter{equation}{0}

When $\ep>0$, all solutions of the problem (\ref{test11}) lie on a unique solution curve joining $(0,0)$ to $(\infty,\infty)$ in $(\mu, w(0))$ plane, see e.g., P. Korman \cite{K2}. When $\epsilon \geq 0.25$, the solution curve is monotone. When $\epsilon$ is small, the solution curve must make at least two turns. Indeed,  this curve begins at $(0,0)$, and it is close to the lower branch of (\ref{test12}) (solutions on the lower branch of (\ref{test12}) are non-singular, hence persist for small $\epsilon >0$), and therefore a turn to the left occurs. After the turn, the solution curve has no place to go for decreasing $\lambda$, so that eventually it must travel to the right, providing us with at least the second turn. Denote by $\beta=\beta(\epsilon)$  the point where $e^{-\frac{1}{\ep + w}}$ changes convexity. When $\ep$ is small, the solution  curve makes exactly one turn to the left in the region where $w(0) \in (0,\beta)$, and it can make a number of turns when $\beta<w(0)<v_0(0)-\epsilon$ (and no turns are possible when $w(0)>v_0(0)-\epsilon$).
\medskip

If the theorem was false, then as $\ep \ra 0$ there would be at least four turning points on every curve. Let $(\mu _1(\ep),w  _1(r,\ep))$ and $(\mu _2(\ep),w  _2(r,\ep))$ denote the second and the third turning points respectively. By Lemma \ref{lma:test10}, both $w  _1(r,\ep))$ and $w  _2(r,\ep))$ are bounded from above by $v_0(0)-\epsilon$ for all $r$. The quantity $w  _1(0,\ep))$ (and hence $w  _2(0,\ep))$) is also bounded from below by, say $\frac14$, for sufficiently small $\epsilon$. Indeed,
\[
\left( e^{-\frac{1}{\ep + w}} \right)''=\frac{e^{-\frac{1}{\ep + w}}}{(w+\epsilon)^4} \left[1-2(w+\epsilon) \right]>0 \,,
\]
for $w<\frac14$ and $\epsilon$ small, and hence only turns to the left are possible when $w(0,\epsilon)<\frac14$.
Take a sequence $\ep _k \ra 0$. Using elliptic estimates, along a subsequence, $(\mu _1(\ep),w  _1(r,\ep))$ tends to a solution $(\eta, v(r))$ of (\ref{test12}). (By the above estimates, $w_1 (r,\epsilon)$ cannot tend to either infinity or zero.) This solution $(\eta, v(r))$  has to be singular, since in any neighborhood of it there are two solutions, that the turning point $(\mu _1(\ep),w  _1(r,\ep))$ brings. (If solution of (\ref{test12}) is non-singular, then by the implicit function theorem there is a unique solution of (\ref{test11}) near it, for $\ep$ small.) 
\medskip

Similarly, $(\mu _2(\ep),w  _2(r,\ep))$ must converge along a subsequence to a singular  solution  of (\ref{test12}), as $\epsilon \rightarrow 0$. Both $(\mu _1(\ep),w  _1(r,\ep))$ and $(\mu _2(\ep),w  _2(r,\ep))$ must converge to the unique turning point of  (\ref{test12}). By Theorem \ref{thm:test10} this  turning point is non-degenerate, which means that for $\ep>0$ small there can be only one singular point of (\ref{test11}) in its neighborhood. But we have two singular points, $(\mu _1(\ep),w  _1(r,\ep))$ and $(\mu _2(\ep),w  _2(r,\ep))$, a contradiction.
\epf
 
\section{Extensions}
\setcounter{equation}{0}

Lemma \ref{lma:test10} does not require $\epsilon$ to be very small, only that $\epsilon <v_0(0) \approx 1.53$. This lemma shows that the set of $w(0)$'s, for which turns may occur for (\ref{test11}), is bounded. It follows that computer generated figures of $S$-shaped bifurcation, like the one in Figure $1$ at $\epsilon=0.22$, can be {\em validated}, i.e., a computer assisted proof of their validity can be given.
\medskip

Other results on perturbation of solution curves can be established similarly.

\begin{thm}
Consider the problem ($u=u(x)$, $x \in R^n$)
\begin{equation}
\label{test20}
\s\s \Delta u+\la (u-\epsilon)(u-b)(c-u) =0 \,, \s \mbox{$|x|<1$}\,, \s  \s u=0 \s \mbox{when $|x|=1$}\,,
\end{equation}
with  constants $0<\epsilon<b<c$, such that $c>2b$. Then for $\epsilon$ sufficiently small, the set of positive solutions of (\ref{test20}) consists of two curves. The lower one starts at $(0,0)$ in the $(\lambda,u(0))$ plane, and it is monotone, tending to $\ep$ as $\la \ra \infty$. The upper curve is parabola-like, with a single turn to the right. Correspondingly, there is a $\la _0 >0$, so that the problem (\ref{test20}) has exactly one positive solution for $\la \in (0,\la _0)$, exactly two  strictly ordered positive solutions at $\la =\la _0$, and exactly three   strictly ordered positive solutions for $\la>\la _0$.
\end{thm}

\pf
By B. Gidas, W.-M. Ni and L. Nirenberg \cite{GNN}, positive solutions of (\ref{test20}) are radially symmetric. Define $f(u)=u(u+\ep-b)(c-\ep-u)$. Then we can write  (\ref{test20}) as 
\begin{equation}
\label{test21}
\s\s \Delta u+\la f(u-\ep) =0 \,, \s \mbox{$|x|<1$}\,, \s  \s u=0 \s \mbox{when $|x|=1$}\,.
\end{equation}
At $\ep=0$, all positive solution of (\ref{test21}) lie on a  parabola-like curve, with a single turn to the right. Moreover, except for the turning point, all solutions are non-singular, and the turning point is non-degenerate, i.e., the conditions of Lemma \ref{lma:test0} hold, and the upper branch is monotone increasing, see \cite{KLO2} and \cite{OS}. As in the proof of the Theorem \ref{thm:*}, we show that this curve preserves its shape  for $\epsilon$ sufficiently small, while monotonicity of the lower solution curve is easy to prove ($f'(u)<0$, when $u<\ep$).
\epf

This is the first exact multiplicity result for a cubic with three positive roots, in dimensions $n>1$. Similar result holds for cubic-like $f(u)$, considered in \cite{KLO} and \cite{OS}.

\end{document}